\documentclass{iopart}
\usepackage{amssymb,setstack}

\newtheorem{Thm}{Theorem}
\newtheorem{Cor}{Corollary}
\newtheorem{Lem}{Lemma}
\newtheorem{Prop}{Proposition}

\eqnobysec
\newcommand{\secref}[1]{\paragraph\ref{#1}}
\newcommand{\lemref}[1]{Lemma~\ref{#1}}
\newcommand{\thmref}[1]{Theorem~\ref{#1}}
\newcommand{\propref}[1]{Proposition~\ref{#1}}

\def\C{{\mathbb C}}
\def\Z{{\mathbb Z}}
\let\paragraph=\S
\def\Mat{\operatorname{Mat}}
\def\qed{\opensquare}

\def\H{{\^H}}
\def\operatorname#1{\rm{#1}}

\def\^#1{\widehat{#1}}
\def\<{\left<}
\def\>{\right>}
\def\bysame{\dash}

\begin{document}
\title{Modules of Abelian integrals and Picard-Fuchs systems}

\begin{abstract}
We give a simple proof of an isomorphism between two
$\mathbb{C}[t]$-modules
corresponding to bivariate  polynomial $H$ with nondegenerate highest homogeneous part:
the  module of relative cohomologies
$\Lambda^2/dH\land \Lambda^1$ and the module of Abelian integrals.
Using this isomorphism, we prove existence and deduce some
properties of the corresponding Picard-Fuchs system.
\end{abstract}

\date{May 27, 2002}
\author{D. Novikov}
\address{Department of Mathematics, Purdue University, West Lafayette, IN 47907,  USA}
\ead{dmitry@math.purdue.edu}
\footnote{The research was supported
by the Killam grant of P. Milman and by James S. McDonnell
Foundation.} \maketitle

Abelian integral is a result of integration of a polynomial
one-form along a cycle lying on level curve (possibly complex) of
a bivariate polynomial considered as a function (possibly
multivalued) of the value of the polynomial. Abelian integrals
appear naturally when studying bifurcations of limit cycles of
planar polynomial vector fields. In particular, zeros of Abelian
integrals are related to limit cycles appearing in polynomial
perturbations of polynomial Hamiltonian vector fields. This is the
reason why sometimes the question about the number of zeroes of
Abelian integrals is sometimes called {\it infinitesimal} Hilbert
16$^{th}$ problem.

The traditional approach to the investigation of Abelian integrals
uses properties of the system of linear ordinary differential
equations satisfied by the Abelian integrals, the so-called
Picard-Fuchs system. This approach is used both in fundamental general finiteness result of
\cite{varchenko,khovanskii} and in exact estimates in the cases of low degree, as in 
\cite{horozov-iliev}. The existence of such a system can be easily
proven due to the very basic properties of branching of Abelian
integrals, see \cite{avg}, and was well known already to Riemann
if not Gauss. Nevertheless an effective  computation of this
system turns out to be a difficult problem. One particular case of
this problem (namely of the hyperelliptic integrals) is quite
classic, see e.g. \cite{roussarie, pham, givental}.  In \cite{redundant} a
generalization of this approach for {\em regular at infinity} (see
below for exact definition) polynomials in two variables is
suggested (in fact, it can be easily generalized for any number of
variables). The main idea of \cite{redundant} is to trade the
minimality of the size of the system (thus redundant) for an
explicitness of the construction and control on the magnitude of the
coefficients.
Another, probably not less important, gain is that the resulting
system is not only Fuchsian, but also has a hypergeometric form.

The control on the magnitude of coefficients in \cite{redundant}
is very important from the infinitesimal Hilbert 16$^{th}$ problem
point of view. Indeed, recent progress towards its solution is
partly based on the principle that solutions of
linear ordinary differential equations
with bounded coefficients cannot oscillate too wildly,
see e.g. \cite{IY:oscill} (simple
proofs of a result of this type can be found in \cite{meandering}
and  \cite{vallee-poussin}). Though more complicated, this
principle still holds for polynomial {\em systems} of differential
equations, see \cite{meandering, lleida} (polynomiality is essential, see
\cite{novikov}). In a slightly modified form, this principle allows
to give results in an upper bound on the number of zeros of an Abelian
integral in terms of the minimal distance between critical values
of its (regular at infinity) Hamiltonian, see \cite{redundant}. As
an application of this principle one can also deduce an effective
upper bound for the number of zeroes
of Abelian integrals corresponding to hyperelliptic Hamiltonians
satisfying some  additional assumption, see \cite{ERA}.

The Picard-Fuchs systems discussed in this paper is irredundant in
the sense that it  has the minimal possible dimension (namely the
dimension equal to the dimension of the homology group
$H^1(\{H(x,y)=t\},\mathbb{C})$ of a  generic fibre). This
minimality allows to guess most of the important information about
the system if the critical values of the Hamiltonian are distinct
and the Hamiltonian is regular at infinity (so-called Morse-plus
Hamiltonians).

We prove existence of such system using decomposition in Petrov
modules. It is easy to see that exact forms  and forms
proportional to $dH$ have zero Abelian integrals, so in fact
Abelian integrals depend on the class of a form in the so-called
{\em Petrov module} --  the quotient of the space of all forms by a
subspace  spanned by exact forms and forms proportional to $dH$, considered in \cite{petrov}.
In \cite{gavrilov} L.Gavrilov proved that the Petrov module of a generic
Hamiltonian is a finitely generated free
$\C[t]$-module.The local counterpart of this statement is
due to E.Brieskorn and M. Sebastiani \cite{brieskorn, sebastiani}.
The proof in \cite{gavrilov} contains a reference to a general nondegeneracy
result (see \cite[Theorem 13.6, Ch III]{avg}),
based on the theory of deformations of Hodge structures.
In the recent preprint  \cite{glutsuk}
the involved constant is computed and exact formulae are given. We
suggest an elementary proof of this result, see also \cite{varchenko1}.

The main idea of \cite{redundant} was to use a connection between
division with remainder of polynomials and differentiation of
Abelian integrals given by Gelfand-Leray formula. In this work we
replace the explicit division with  remainder by decomposition in
Petrov modules in order to get the same result. This is still
enough for the construction of the system, though the result is
less explicit. Yet one can still guess all singular points and get
some information about coefficients. However, we show that the
resulting irredundant system is not always Fuchsian, namely it can
have regular but non-Fuchsian point at infinity. Though after a suitable
rational gauge transform the irredundant system becomes  Fuchsian
(see \secref{fuchs}), the nice form of
\thmref{irred} is then lost.

\subsection{Acknowledgments}
I am grateful to S. Yakovenko for the numerous discussions and
help in preparation of this text. I am grateful to Yulij
Ilyashenko for an opportunity to read his recent preprint and
numerous discussions. I am grateful to anonymous  referee for
several important remarks.

\section{Genericity and generalities}
In what follows we always assume that our polynomial $H(x,y)$ is
regular at infinity, i.e., that its highest homogeneous part
$\H(x,y)$ is a product of {\it pairwise different} linear factors.

One can easily prove that for regular at infinity polynomial
$H(x,y)$ of degree $n+1$  the homogeneous polynomial $\H$ has an
isolated critical point (necessarily of multiplicity $\mu=n^2$) at
the origin $(x,y)=(0,0)$, its level curves $\{\H=c\}\subset\C^2$
are nonsingular for $c\not=0$. Moreover, the level curves of $H$
intersect transversally the line at infinity, and foliation of
$\C^2$ by level curves of $H$ is locally topologically trivial
over $\C\setminus\Sigma$, where $\Sigma$ is the set of $\le(\deg
H-1)^2$ critical values of $H$. In other words, the only atypical
values are the critical ones.

By Abelian integral we mean a result  of integration of a one-form
$\omega$ along a continuous family of cycles
$\delta(t)\subset\{H=t\}$ considered as a function of $t$:
$$
I_{\omega,\delta}(t)=\oint_{\delta(t)}\omega.
$$

Basic properties of Abelian integrals can be found in \cite{avg}.
We will need the following ones. First,  Abelian integral depends
not on $\delta$ itself but on its class of homology $[\delta]\in
H_1(\{H=t\},\Z)$ only. Also, the Abelian integral corresponding to
a form $\omega$ is identically zero if $\omega=fdH+dg$, i.e.
Abelian integrals depend in fact on the relative cohomology class of
$[d\omega]\in\Lambda^1/dH\land \Lambda^0+d\Lambda^0$ only. This quotient
module -- the so-called Petrov module -- is a $\mathbb{C}[t]$-module with respect to a standard
multiplication $t[\omega]=[H(x,y)\omega]$.

Second, the Abelian integrals are holomorphic multivalued
functions of a complex variable $t$ branching at the critical
values of $H$ only (for $H$ regular at infinity). So the space of
all Abelian integrals is also a $\mathbb{C}[t]$-module with
respect to a natural multiplication by $t$.

We will prove that these two modules coincide for a regular at infinity
polynomial $H(x,y)$. We prove existence  of the corresponding
Picard-Fuchs system using this isomorphism, and find out some of
its properties.

\section{Nondegeneracy of the principal determinant}\label{nondegeneracy}
Here we prove that the homogeneous forms generating the
$\mathbb{C}[t]$-module $\Lambda^2/dH\land \Lambda^1$, also
generate the first cohomology group of a generic level curve of
$H$. Note that this $\mathbb{C}[t]$-module is isomorphic to
$\mathbb{C}[x,y]/\<H_x,H_y\>$. Indeed, $Qdx\land dy=Rdx\land
dy+(H_x dx+H_ydy)\land (Adx+Bdy)$ is equivalent to
$Q=R+BH_x-AH_y$.

Recall that nonsingular level curves of a regular at infinity polynomial $H$
of degree $n+1$ carry $\mu=n^2$ vanishing cycles $\delta_j(t)$
that generate the whole first homology group of all regular curves
$\{H=t\}$ \cite{avg,gavrilov}. For any collection of $\mu$
polynomial 1-forms $\omega_1,\dots,\omega_\mu$ the {\em period
matrix\/} $X_\omega(t)$ formed by integrals of $\omega_i$ over
$\delta_j$ (integrals of the same form occur in the same row, the
same cycle corresponds to entries of the same column) has the same
monodromy. The monodromy transformations act on $X(t)$ as
multiplications from the right by constant monodromy matrices that
are unimodular by virtue of Picard--Lefschetz formulas
\cite{avg}. Thus $\det X(t)$ is a single-valued function that must
have zeros at all critical values $t=t_j$, $j=1,\dots,\mu$, counting multiplicities,
since the columns corresponding to the cycles vanishing at $t_j$ become
zero at $t_j$ (we use the fact that to any critical value of multiplicity $\nu$ correspond $\nu$ linearly independent
cycles vanishing at this critical value). 
As the growth of $X$ at infinity is at most
polynomial, $\det X(t)$ is a polynomial divisible by
$\Delta_H(t)=\prod_{j=1}^\mu (t-t_j)$.

\begin{Lem}[cf.~with \cite{gavrilov}, Lemma~2.2]\label{nondeg}
If the $2$-forms $d\omega_i$ generate $\Lambda^2/dH\land
\Lambda^1$ and
$$
 \sum_{i=1}^\mu \deg \omega_i=\mu\deg H,
$$
then $\det X_\omega(t)=c(t-t_1)\dots(t-t_\mu)$ with $c\ne 0$ (some $t_i$ may coincide).
\end{Lem}
The constant $c$ depends both on the choice of $\omega_i$ and on
the choice of the cycles $\delta_j(t)$. Its actual calculation is
a difficult task, see \cite{glutsuk}.

\noindent{\bf Remark.} The condition on the degrees of the forms
in the \lemref{nondeg} is automatically satisfied if the
Hamiltonian $H(x,y)$ and the forms $d\omega_i$ are homogeneous,
see \cite{avg}. For non-homogeneous Hamiltonian this condition is essential. 
Among other things, this condition implies that the
highest homogeneous parts $\^\omega_i$ of $d\omega_i$ form a basis
of $\Lambda^2/d\H\land \Lambda^1$. Vice versa, any monomial basis of 
$\Lambda^2/d\H\land \Lambda^1$ is a basis of $\Lambda^2/dH\land
\Lambda^1$ satisfying this condition (and this is a standard way to get 
a basis of $\Lambda^2/dH\land \Lambda^1$).

The proof is based on the calculation of the ``principal term'' of
the asymptotic of $X(t)$ at infinity.

\begin{Lem}\label{tDC}
For any collection of  polynomial 1-forms $\omega_i$ the period
matrix $X_\omega(t)$ admits a converging expansion
\begin{equation}\label{exp-inf}
  X_\omega(t)=t^DC(t),
  \qquad C(t)=\sum_{k=0}^\infty C_k t^{-k/(n+1)},
\end{equation}
where $D$ is the diagonal matrix with the entries $d_i=\deg
\omega_i/(n+1)$, $C_0,C_1,\dots$ are constant matrices and
$C_0=C(\infty)$ is the matrix of integrals of the highest
homogeneous parts $\^\omega_i$ of forms $\omega_i$ over vanishing
cycles lying on the level curve $\{\H=1\}$.
\end{Lem}

\noindent{\em Proof of \lemref{tDC}\/.} The level curve
$\{H(x,y)=t\}$ in the variables $x=t^{1/(n+1)}\^x$,
$y=t^{1/(n+1)}\^y$ becomes a family of the curves
$$
 \H(\^x,\^y)+t^{-1/(n+1)}H_{n}(\^x,\^y)+\cdots=1,
$$
where the left hand side is a polynomial in $\^x,\^y$ and
$t^{-1/(n+1)}$. In other words, we have an analytic in
$t^{-1/(n+1)}$ perturbation of the limit curve
$\{\H(\^x,\^y)=1\}\subset\C^2$ that is nonsingular (since $\H$ has
no multiple factors). Integrals of any (constant or analytic in
$t^{-1/(n+1)}$) 1-form over any continuous family of cycles on
such family will be also analytic in $t^{-1/(n+1)}$.

The forms $\omega_i$ after rescaling become
$t^{d_i}(\theta_i+t^{-1/(n+1)}\eta_i)$, where $d_i={\deg
\omega_i/(n+1)}$, $\theta_i$ is a new independent of $t$
homogeneous polynomial form (corresponding to the highest
homogeneous part $\^\omega_i$ of $\omega_i$) and $\eta_i$ is
another polynomial form. Therefore the integrals of $\omega_i$
over cycles $\delta_j(t)$ on the level curves $\{H=t\}$ can be
expanded in the converging series in $t^{-1/(n+1)}$ of the form
$$
 \oint_{\delta_j(t)}\omega_i=
 t^{d_i}(c_{0,ij}+c_{1,ij}t^{-1/(n+1)}+\cdots),
$$
if  $c_{0,ij}$ is the integral of $\^\omega_i$ over the cycle
$\delta_j\subset\{\H=1\}$. \qed

\noindent{\bf Remark.} The representation (\ref{exp-inf}) is unique
only if we fix the diagonal matrix $D$. Otherwise the power $t^D$
may itself be expanded as a series in powers of $t^{-1/(n+1)}$,
thus yielding an essentially different representation.

\begin{Cor}
The determinant of the period matrix $X_\omega(t)$ is a polynomial
of degree at most $m=m(\omega)=\tr D=\sum_i \deg\omega_i/(n+1)$.
If this number is not integer, then automatically $\det C_0=0$ for
this choice of the forms, otherwise the leading term $t^m$ of
$\operatorname{det} X_\omega(t)$ enters with the coefficient $\det
C_0$.
\end{Cor}

\noindent{\em Proof of the \lemref{nondeg}\/}. Given the
assumption on the degrees $\deg \omega_i$, the determinant $\det
X_\omega(t)$ is a polynomial of degree $\le \mu$, and hence (by
the divisibility property noted above) it must have a form
$c\prod(t-t_j)$. We need only to verify that $c\ne0$, and from the
asymptotic formulas we see that $c=\det C_0$, so our goal is to
prove that $C_0$ is a nondegenerate matrix.

The calculation above shows that the matrix $\^X(t)$ of periods of
$\^\omega_i$ over the level curves of a {\em homogeneous\/} part
$\H$, can be represented as $t^DC_0$ (the same expansion without
inferior terms). Thus if $C_0$ is degenerate, then there exists a
linear combination $\^\delta(t)=\sum_1^\mu r_j\delta_j(t)$,
$r_j\in\C$, of vanishing cycles on the level curves of $\H$, such
that integrals of all forms $\^\omega_i$ over the cycle
$\^\delta(t)$ are identically zeros.

Take any polynomial 2-form $d\omega$. Since the forms
$d\^\omega_i$ form a basis of $\Lambda^2/d\H\land \Lambda^1$, the
form $d\omega$ can be divided out by $d\^H$ with remainder in the
span of $d\^\omega_i$, i.e.,
$$
 d\omega=d\H\land \eta+\sum_1^\mu c_id\^\omega_i,\qquad
 c_i\in\C,
$$
where $\eta$ is a suitable polynomial 1-form.

This representation is not unique. However, since $H$ is regular
at infinity, one can construct such representation with  degree of
$\eta$ being less than $\deg\omega$ (in fact, less or equal to
$\deg\omega-\deg H$, see \cite{redundant}).

Recall that the derivative of an Abelian integral of a form
$\omega$ with respect to $t$ is again an Abelian integral of the
{\em Gelfand-Leray residue} $\theta=\frac{d\omega}{dH}$ of the
form  $\omega$:
$$
\frac d {dt} \oint_{\delta(t)}\omega=\oint_{\delta(t)}\theta,
$$
if $dH\land \theta=d\omega$.

Return to the division with remainder of the form $d\omega$ by
$d\H$.  Integrating over the cycle $\^\delta(t)$ and using the
Gelfand--Leray formula, we see that
$$
 \frac  d{dt}\oint_{\^\delta(t)}\omega=\oint_{\^\delta(t)}\eta,
 $$
since integrals of $\frac{d\^\omega_i}{d\H}$ over $\^\delta(t)$
all vanish. In other words, the derivative of any Abelian integral
of a polynomial form over the cycle $\^\delta(t)$ is again an
Abelian integral of a {\em polynomial\/} form. Since the cycle
$\^\delta(t)$ is also vanishing at $t=0$ (recall that we deal with
the homogeneous case and all $\delta_i(t)$ vanish at the same
value $t=0$), the limit of $\oint_{\^\delta(t)}\eta$ is zero for
any polynomial form $\eta$ as $t\to0$.

As the Gelfand--Leray derivative $\eta$ is a polynomial form of
smaller degree, the above argument can be repeated, showing that
some derivative of the initial integral
$\oint_{\^\delta(t)}\omega$ is zero. Since the integral itself and
all its derivatives tend to zero as $t\to0$, we conclude that the
initial Abelian integral is identically zero. Since $\omega$ was
arbitrary, this proves that {\em integrals of all polynomial forms
over the cycle $\^\delta(t)$ are identically zeros\/}.

But this is clearly impossible unless $\^\delta\equiv 0$ in $H_1(\{\H=1\},\mathbb C)$.
The shortest way to show this is to refer to \cite{avg}, where the following
statement is proved.

\begin{Lem}[\cite{avg}]
For an isolated singularity with Milnor number $\mu$ one can
always construct $\mu$ holomorphic 1-forms
$\theta_1,\dots,\theta_\mu$ such that the period matrix
$X_\theta(t)$ (integrals of $\omega_i$ over all vanishing cycles)
will have the determinant $\det X_\theta(t)=t^\mu+\cdots$. \qed
\end{Lem}

This lemma can be applied to the homogeneous germ $\H$ and the
forms in \cite{avg} are in fact constructed polynomial (of course,
of sufficiently high degrees). Namely, for an arbitrary
nonzero cycle (in particular, for $\^\delta(t)$) a suitable
linear combination of $\theta_i$ has integral not identically
zero, which contradicts the choice of $\^\delta(t)$. \qed

\noindent{\bf Remark.} The assertion of the above Lemma is by far
much stronger than required to complete the proof: it would be
sufficient to find just one polynomial  form in $\mathbb C^2$ such that the integral of
its restriction to  the affine curve $\{\H=t\}$ along $\^\delta(t)\ne0$ would be  non-zero.
This can be done using the fact that $\{\H=1\}$ is a Stein
manifold, and therefore each element of its cohomology group can
be realized as a restriction of a holomorphic one-form on $\C^2$.
More exact, let $\omega$  be a holomorphic form on $\{\H=1\}$ such that its integral along
$\^\delta $ is nonzero. One can find a holomorphic form $\tilde\omega$
on $\mathbb C^2$ which restriction to $\{\H=1\}$
is cohomologous to $\omega$. Since the cycles generating $H_1(\{\H=1\}, \C)$ have  compact
representatives,
a  polynomial
one-form sufficiently close approximating  $\tilde\omega$  on  a sufficiently big compact
 will also produce nonzero integral along
$\^\delta(t)$  (since  analytic in $\mathbb C^2$
coefficients of the form $\tilde\omega$  can be uniformly approximated by
polynomials on any compact set).

\section{Module of the Abelian integrals}
 Now, after \lemref{nondeg} is proved,
we can immediately prove that integrals of the forms $d\omega_i$
generate over $\C[t]$ the entire module of Abelian integrals. The
proof appears in \cite{gavrilov} and is a
straightforward application of the Cramer rule. We reproduce this
proof here for reader's convenience.

\begin{Prop}[Gavrilov theorem \cite{gavrilov}]\label{module}
Let  $\omega_1,\dots,\omega_\mu$ be one-forms such
that $\sum_{i=1}^{\mu}\deg\omega_i=\mu\deg H$, and suppose that
the polynomials $\frac{d\omega_i}{dx\land dy}$ are linearly
independent modulo the gradient ideal $<H_x,H_y>$ in $\C[x,y]$.

Then integral of any polynomial $1$-form $\omega$ can be
represented as a linear combination of integrals of the forms
$\omega_i$ with polynomial in $t$ coefficients\/{\rm:} for any
cycle $\delta(t)$ on the level curve $\{H=t\}$
\begin{equation}\label{petrov-module}
\fl \oint_{\delta(t)}\omega=\sum_{i=1}^\mu
  p_i(t)\oint_{\delta(t)}\omega_i,\quad p_i(t)\in\C[t],\ (n+1)\deg p_i+\deg \omega_i\le \deg\omega.
\end{equation}

\end{Prop}

\noindent{\bf Remark.} The condition on degrees is again essential: if
$H(x,y)$ is not homogeneous, then not every basis of monomial forms
of $\Lambda^2/dH\land\Lambda^1$ generates the
Petrov module.  A (more transparent
weight-homogeneous) example is $H=y^2+x^4-x^2$ and the set of
monomial forms $dx\land dy, x^2 dx\land dy, x^5 dx\land dy$.
However,for homogeneous $H$ and homogeneous $\omega_i$ this condition is
satisfied automatically, see \cite{avg}.

\noindent{\em Proof\/.} We look for a tuple of real functions
$p_i(t)$ such that identically over $t$ and for any vanishing cycle
$\delta(t)=\delta_j(t)$ the equality (\ref{petrov-module})
holds. These equations for each $t$ form a linear nonhomogeneous
system with the matrix of coefficients $X(t)$ being the period
matrix $\oint_{\delta_j}\omega_i$ and the column of right hand
sides being periods of the form $\omega$.

Since the matrix $X(t)$ is nondegenerate (for all $t\ne t_j$), the
solution of this system can be found by the Cramer rule: each
$p_i$ is a ratio of two determinants. The denominator is $\det
X(t)=c\prod_j(t-t_j)$, whereas the numerator is the determinant of
the period matrix obtained by replacing $\omega_i$ by $\omega$. By
the same arguments as in the beginning of \secref{nondegeneracy},
the numerator should be a polynomial divisible by
$\prod_1^\mu(t-t_j)$, hence the inequality $c\ne0$ ensures that
the ratio is in fact a polynomial function of $t$. To estimate the
degree of the nominator, we use Corollary to \lemref{tDC}:  it is
no greater than $\deg\det
X(t)+\frac{\deg\omega-\deg\omega_i}{n+1}$. Therefore $\deg p_i\le
(\deg\omega-\deg\omega_i)/(n+1)$.  \qed

\noindent{\bf Remark.} The uniqueness of the representation
(\ref{petrov-module}) follows from a theorem by Gavrilov (see
\cite{gavrilov,gavrilov-2}) that a polynomial 1-form with all zero
periods must be necessary $a(x,y)dH+db(x,y)$, where $a,b$
appropriate polynomials, provided that the Hamiltonian $H(x,y)$ is regular
at infinity (the conditions in \cite{gavrilov} are
even weaker). This result is a generalization of an earlier  result of
Ilyashenko \cite{ilyashenko1}.

The local counterpart of \propref{module} claims that the ring of
relative cohomology is finitely generated as a $C\{t\}$-module
(Brieskorn--Sebastiani \cite{brieskorn,sebastiani}).

\section{Derivation of the irredundant Picard--Fuchs system and its
elementary properties} \label{system}

Let $\omega_1,\dots,\omega_\mu$ be
polynomial 1-forms as in \propref{module}, i.e., they satisfy the condition $\sum_1^\mu
\deg\omega_i=\mu\deg H$ and their differentials $d\omega_i$ generate
$\Lambda^2/dH\land \Lambda^1$.

The second assumption guarantees that
we may divide out the 2-forms
$H(x,y)d\omega_i$ for all $i=1,\dots,\mu$, obtaining
\begin{equation}\label{division}
 H\,d\omega_i=dH\land \eta_i+\sum_{j=1}^\mu a_{ij}d\omega_j,\qquad
 i=1,\dots,\mu,
\end{equation}
with appropriate polynomial forms $\eta_i$ of degrees
$\deg\eta_i\le \deg\omega_i\le 2n$. This by the Gelfand--Leray
formula implies that for any cycle $\delta(t)$
$$
 (t-A)\dot I(t)=J(t), \quad{\rm where}
\quad I=(\oint_{\delta(t)}\omega_1,...,\oint_{\delta(t)}\omega_{\mu})^T,
\quad J=(\oint_{\delta(t)} \eta_1,\dots,\oint_{\delta(t)} \eta_\mu)^T.
$$

Here occurs the difference with the computations from
\cite{redundant}: we cannot claim that the integrals $J_i$ are
linear combinations of $I_j$, since the {\em linear span\/} of the
forms $d\omega_i$ does not contain  all $2$-forms of
degrees $\le 2n$ (in \cite{redundant} this decomposition was written for {\em
all\/} monomials of degree $\le 2n$ which resulted in a hypergeometric system
of doubled size with a Fuchsian singularity at infinity).

However we can use the decomposition provided by \propref{module}
and write
$$
 J(t)=B(t)I(t),\qquad B(t)=B_0+tB_1,
$$
i.e  $B(t)$ is a  matrix polynomial of degree $\le 1$.

This proves the following result.

\begin{Thm}\label{irred}
The period matrix $X(t)$ of the forms $\omega_i$ satisfying the
above three conditions, is a nondegenerate solution to the system
of first order linear ordinary differential equations
\begin{equation}\label{irred-sys}
  (t-A)\dot X(t)=(B_0+B_1t)X(t),\qquad
  A,B_0,B_1\in\Mat_{\mu\times\mu}(\C).
\end{equation}
\end{Thm}

Some properties of the matrices $A,B_0,B_1$ can be established by
a simple inspection. First, after  identification of $\Lambda^2/dH\land \Lambda^1$
with $\C[x,y]/<H_x, H_y>$, the equation (\ref{division})
means that $A$ is a matrix of multiplication by
$H$ in  $\C[x,y]/<H_x, H_y>$.

Suppose for a moment that $H(x,y)$ has $\mu$ simple pairwise different critical values.
Let $(x_j,y_j)$, $j=1,\dots,\mu$ be critical
points of $H$. Denote by $\mathbf v_j$ the $\mu$-dimensional
vector, whose components are the coefficients
$\frac{d\omega_i}{dx\land dy}$ evaluated at the point $(x_j,y_j)$.
Such vectors form a basis in $\C^\mu$ by the second condition
imposed on the forms. For example, if the coefficients of $\omega_i$ 
are monomials $x^\alpha y^\beta$
with $0\le \alpha,\beta\le n-1$, then together $\mathbf v_j$,
$j=1,\dots,\mu$ form a two-dimensional analog of the Vandermonde
matrix.

\begin{Prop}\label{A}
The matrix $A$ is diagonalizable, its eigenvalues are critical
values of $H$ whereas the eigenvector corresponding to the
critical value $t_j$ is $\mathbf v_j$.
\end{Prop}

\noindent{\em Proof\/.} The right hand side of the expression
(\ref{irred-sys}) has $j$-th column zero if evaluated at the point
$t=t_j$, since the corresponding cycle vanishes. The corresponding
column of the matrix $\dot X(t_j)$ is therefore in the kernel of
$(t_j-A)$. Since the number of critical values is equal to the
dimension of the system (recall we are dealing with the
irredundant case), this proves the assertion about
diagonalizability and the spectrum of $A$.

To complete the proof we need only to compute the derivatives
$\dot I_i(t_j)$. The Gelfand--Leray derivative $d\omega_i/dH$ has
zero residues on all nonsingular level curves, but restricted on
$\{H=t_j\}$ it has a nontrivial residue. This can be immediately
seen for the normal form when $H(x,y)=y^2-x^2$ (note that all
considerations are local, so one can use the Morse normal form
near the critical point $(x_j,y_j)$). Indeed, if
$d\omega=f(x,y)\,dx\land dy$, then $d\omega/dH$ can be chosen as
$\frac12f\,dx/y$, and its restriction on (one of the two smooth
branches of) the curve $H=0$, say, $y=x$, yields a meromorphic
1-form $\frac12 f(x,x)\,dx/x$, whose residue (integral over a
small loop around $x=0$) is $\pi if(0,0)$. Returning to the
initial problem, we see that $\oint_{H=t_j}\frac{d\omega_i}{dH}$
differs from $\pi\frac{d\omega_i}{dx\land dy}{(x_j,y_j)}$ by a
nonzero factor, the Hessian of the transformation taking $H$ into
the Morse form as above. Since this nonzero factor is common for
all forms, we see that the vector of residues $(\dot
I_1(t_j),\dots,\dot I_\mu(t_j))$ is proportional to the vector
$\mathbf v_j$ whose coordinates are $\frac{d\omega_i}{dx\land
dy}{(x_j,y_j)}$, $i=1,\dots,\mu$. \qed

By continuity one can conclude that 
\begin{Cor}
For any regular at infinity Hamiltonian $H(x,y)$   its  critical
values $t_j$ counted with multiplicities are the eigenvalues of 
the matrix $A$, and the vectors $\mathbf v_j$ are eigenvectors of $A$
\end{Cor}

The matrices $B_0$, $B_1$ in principle can be computed by
evaluating the expansion for $X(t)$ at infinity, see \lemref{tDC}.
One can guess some of their properties just by taking $d\omega_i$
homogeneous and of nondecreasing degree.
\begin{Prop}\label{B}
Let $d\omega_i$ be homogeneous and $\deg d\omega_i\le\deg
d\omega_j$ whenever $1\le i<j\le\mu$. Then $B_0$ and $B_1$ are
both lower triangular. Moreover, the diagonal entries of $B_0$ are
just the degrees of the forms divided by $\deg H$, and $B_1^2=0$
\end{Prop}

\noindent{\em Proof\/.} This follows from the careful analysis of
the forms $\eta_i$ in (\ref{irred-sys}). Indeed, $\deg\eta_i\le
\deg\omega_i$, so in the decomposition of $\eta_i$ provided by
\propref{module} appear only forms of degree not greater than
$d\omega_i$. Moreover, it is easy to see (using Euler identity)
that the highest homogeneous term of $\eta_i$ is equivalent in the
Petrov module to $\frac {\deg\omega_i}{\deg H}d\omega_i$, see
\cite{redundant}. This together implies that $B_0$ is lower
triangular with prescribed diagonal elements. From the same
estimates of the \propref{module} follows that entries
$(B_1)_{ij}$ of the matrix $B_1$ can be nonzero only if
$\deg\omega_i-\deg\omega_j\ge\deg H$, so $B_1$ is lower triangular
and, since $\max_{ij}(\deg\omega_i-\deg\omega_j)=2\deg H -4<2\deg
H$, already $B_1^2=0$.

\begin{Cor}
The matrix $B_0+tB_1$ is invertible for all $t$.
\end{Cor}

\section{ Picard-Fuchs system can be non-Fuchsian}\label{fuchs}
From the analysis above follows that all finite singular points of
the system (\ref{irred-sys}) coincide with the critical values of
$H$. Moreover, all finite singularities turn out to be Fuchsian
for Morse-plus $H(x,y)$ (which, by definition, means that the
matrix $(t-A)^{-1}(B_0+tB_1)$ of coefficients of the system of the
\thmref{irred} has poles of the first order). Indeed, the
Picard-Fuchs system has a Fuchsian singularity at $\lambda_i$ if
and only if the matrix $(t-A)^{-1}$ has a simple pole at
$\lambda_i$ (due to invertibility of $B_0+tB_1$ for all $t$). This
is equivalent to the diagonalizability of the matrix $A$,  so is
true for Morse-plus Hamiltonian $H(x,y)$.

For a general regular at infinity $H(x,y)$ the finite singular points can
be non-Fuchsian. Indeed, the matrix $A$ is the matrix of
multiplication by $f$ in $\C[x,y]/<H_x, H_y>$, and this ring is a
direct sum over all critical points of $H(x,y)$ of the
corresponding local rings (\cite[Max Noether's $AF+BG$ Theorem]{GH}). It
follows that  $A$ is diagonalizable if and only if the operator of
multiplication by $H$ is diagonalizable in each local ring. This
is true if and only if the germ of $H(x,y)$ is (equivalent to)
quasi-homogeneous, see \cite{avg}, so fails in general.

Also, unless $B_1=0$, the singular point at infinity is
non-Fuchsian. This is also possible, see below an example.

It is easy to see that the Picard-Fuchs system of the
\thmref{system} is {\em equivalent} to a Fuchsian one. Indeed, due
to the regularity at infinity assumption the monodromy of the
irredundant system corresponding to a circle around infinity is
diagonalizable, so this equivalence is a particular case of a
positive solution (essentially due to Plemelj) of the Riemann-Hilbert problem
in the case of diagonalizability of  one of the local monodromies,
see \cite{dynsyst1}. Moreover, in \cite{kostov} it is proved,
modulo a conjecture due to Bolibruch, that this system is
equivalent to a Fuchsian one for any $H(x,y)$, even degenerate
ones. However, the equivalent system will not have the fairly
simple form of \thmref{system}.

Here is an example of a Hamiltonian with a nonzero matrix $B_1$.

\noindent{\bf Example.} Consider the Hamiltonian
$H(x,y)=x^5+y^5+x^2y^2+ax+by$. For a suitable choice of $a,b$ this
Hamiltonian is Morse-plus. As a basis of the quotient
$\Lambda^2/dH\land\Lambda^1$ we  take the forms
$d\omega_{ij}=x^iy^jdx\land dy$, $0\le i,j\le 3$. We will show
that any form $\eta$ defined by the
decomposition
$$
Hd\omega_{33}=dH\land\eta+\sum_{0\le i,j\le
3}a_{ij}d\omega_{ij}
$$
is equivalent  to $\frac 1 {175} t\omega_{00}+\sum_{0\le
i,j\le3}\beta_{ij}\omega_{ij}$  in the  Petrov module corresponding to $H(x,y)$,
with $\beta_{ij}$ being constant (so the matrix
$B_1$ has a nonzero entry equal to $\frac 1 {175}$).

Although $\eta$ is defined non-uniquely by the Gelfand-Leray
formula above, its Abelian integrals do (and therefore its class in the Petrov module).
So we can use any $\eta$
we like. Applying the ``division with remainder'' algorithm of
\cite{redundant} we find the first terms of a form $\eta$ solving
the equation above:
\begin{eqnarray}
\eta=\frac 1 5 (x^3y^3(xdy-ydx)+
(\frac 1 {175}xy^5dy-\frac 6 {175}x^5ydx)+\eta_1=\nonumber\\
=\frac {x^3y^3} 5 (xdy-ydx) + \frac {y^5+x^5} {175}
xdy-\frac{d(x^6y)}{175}+\eta_1,\nonumber
\end{eqnarray}
where by $\eta_1$ we denote forms of  degree less than $7$.  It is
easy to see that in the Petrov module  the
first term is
 equivalent to $\frac 8 5 \omega_{33}$ and the second term is equivalent
to $\frac 1 {175} Hxdy+ \frac 1 {175}(x^5+y^5-H)xdy=t\frac 1 {175}
xdy+\eta_2$. Since the degrees of both $\eta_1$ and $\eta_2$ are
less than $7$, the form $\eta_1+\eta_2$  is equivalent in the
$\C[t]$-module of Abelian integrals to a linear combination with
constant coefficients of forms $\omega_{ij}$, by virtue of
estimates of the Corollary to \lemref{tDC}.\qed

\end{document}